# Stratified Picard–Lefschetz theory

Victor A. Vassiliev*

To Robert MacPherson

**Abstract**

The monodromy action in the homology of level sets of Morse functions on stratified singular analytic varieties is studied. The local variation operators in both the standard and the intersection homology groups defined by the loops around the critical values of such functions are reduced to similar operators in the homology groups of the transversal slices of the corresponding strata.

## Introduction

The classical Picard–Lefschetz theory deals with the problems of the following kind. Let $A \subset \mathbf{CP}^n$ be a smooth analytic variety and $\mathrm{reg}(A)$ the set of all hyperplanes $X \in (\mathbf{CP}^n)^*$ transversal to $A$. For all $X \in \mathrm{reg}(A)$ the spaces $A \cap X$ have the same topological type, and the group $\pi_1(\mathrm{reg}(A))$ acts naturally on the space $H_*(A \cap X)$.

Many problems of integral geometry and partial differential equations can be reduced to the study of this monodromy action or of the similar action on the groups $H_*(X \setminus A)$, $H_*(\mathbf{CP}^n \setminus (A \cup X))$, $H_*(\mathbf{CP}^n, (A \cup X))$ etc.: this action is responsible for the ramification of many special functions and fundamental solutions given by integral representations, see f.i. [Pham], [ABG], [AVGL], [AV], [V].

The group $\pi_1(\mathrm{reg}(A))$ is generated by finitely many *simple loops* embracing the set $(\mathbf{CP}^n)^* \setminus \mathrm{reg}(A)$ close to the simplest points of it (which correspond to planes simply tangent to $A$). The monodromy action of any such loop is described by the classical Picard–Lefschetz formula.

In this article we study similar problems in the case when the variety $A$ is singular. The main result is the stabilization theorem (subdivided into Theorems 1–10 below) expressing the ramification of groups $H_*(A \cap X)$ and $IH_*(A \cap X)$ close to planes tangent to singular strata of $A$ in the terms of transversal slices of these strata.

Let $\sigma$ be a $k$-dimensional stratum of some analytical Whitney stratification of $A$, and $l$ a small loop in $\mathrm{reg}(A)$ close to a plane $Y \in (\mathbf{CP}^n)^* \setminus \mathrm{reg}(A)$ simply tangent to $\sigma$. We reduce the monodromy action of the loop $l$ on the group $H_*(A \cap X)$ to a similar action on the group $H_*(A \cap X \cap T)$ where $T$ is an $(n-k)$-dimensional plane transversal to $\sigma$, see Figure 1. In particular, for smooth $A$ we get a new proof of the Picard–Lefschetz formula. Our theorem is thus in the same relation with the standard Picard–Lefschetz theorem as the main theorem of the first part of [GM 86] with the standard Morse theory.

*Mathematics College of Independent Moscow University and Institute for System Studies, Moscow, Russia. e.mail merx@glas.apc.org (Subject: for V. A. Vassiliev)  Research supported by International Science Foundation Grant MQO 000  *Submitted to Selecta Mathematica*

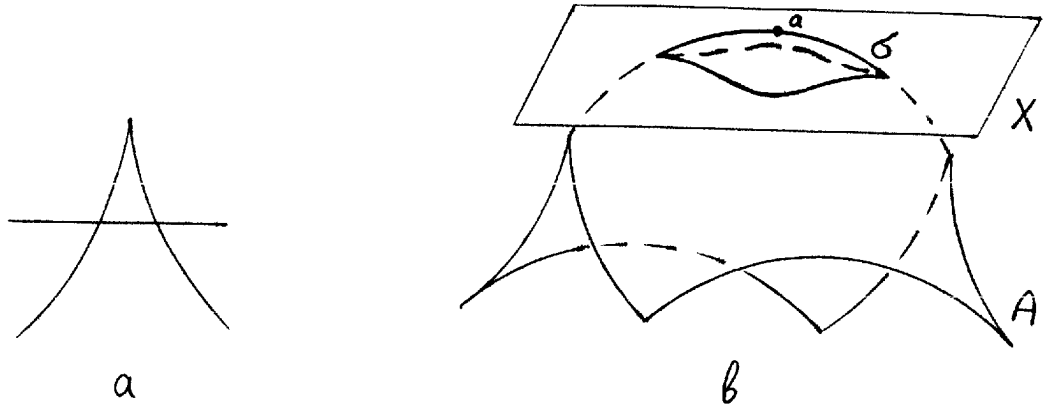

Figure 1: A Morse function on a stratified variety and transversal slice

We also prove similar stabilization formulae for the intersection homology group (of middle perversity) of $A \cap X$ and for the intersection form on it.

This work was motivated by the Newton–Arnold problem on integrability of convex bodies, see [Newton], [Arnold], [AV], and especially by the search of geometrical obstructions to the integrability. For some results obtained in this way see [V], [AVGL 2].

**On the notation.** We consider only the homology groups reduced modulo a point in the case of absolute homology and modulo the fundamental cycle in the case of relative homology of a complex analytic variety; this reduction is emphasized by the tilde: $\tilde{H}_*$. We often use a short notation of type $\tilde{H}_*(X, Y)$ instead of a more rigorous $\tilde{H}_*(X, Y \cap X)$ or $\tilde{H}_*(X \cap Y, Y)$. The sign □ denotes the end or absence of a proof.

# 1 Stating the problem

Let $A$ be a complex analytic subvariety in $\mathbf{C}^n$ with a fixed analytic Whitney stratification (see f.i. [GM 86]), let $\sigma \subset A$ be a stratum of dimension $k$, $a$ a point of $\sigma$, and $B \equiv B_\varepsilon$ a small closed disk in $\mathbf{C}^n$ centred at $a$, so that $B$ has nonempty intersections only with the strata adjoining $\sigma$ (or $\sigma$ itself) and $\partial B$ is transversal to the stratification. In particular, the induced partition of the pair $(B, \partial B)$ is again a Whitney stratification of $B$.

Let $\tilde{f} : (B, a) \to (\mathbf{C}, 0)$ be a holomorphic function such that $d\tilde{f} \neq 0$ in $B$. Denote by $f$ the restriction of $\tilde{f}$ on $A$.

**Definition** (cf. [GM 86]). The function $f$ has a *Morse singularity* at $a$ if its restriction on the manifold $\sigma$ has a Morse singularity at $a$ and for any stratum $\tau \neq \sigma$ and any sequence of points $b_i \in \tau$ converging to $a$ and such that the planes tangent to $\tau$ at $b_i$ (considered as the points in the associated Grassmann bundle $G^{\dim \tau}(T_*\mathbf{C}^n)$) converge to a plane in $T_a\mathbf{C}^n$, this limit plane does not lie in the hyperplane $\{d\tilde{f}|_a = 0\} \subset T_a\mathbf{C}^n$.

Let $f$ be such a function. For any $\lambda \in \mathbf{C}^1$ consider the set $X_\lambda \equiv \tilde{f}^{-1}(\lambda) \cap B$. If $B$ is sufficiently small and $|\lambda|$ is sufficiently small with respect to $B$, then $X_\lambda$ is transversal to all strata of $A \cap B$ outside the point $a \in X_0$, in particular $a$ is the unique critical point of $f$ in $B$. By the Thom's isotopy theorem, for a very small neighbourhood $D \equiv D_\delta$ of the point $0 \in \mathbf{C}^1$ the projection $f : A \cap B \to \mathbf{C}^1$ defines a locally trivial fibre bundle over $D \setminus 0$ with fibres $A \cap X_\lambda$, and the restriction of this projection on $B_\varepsilon \setminus B_{\varepsilon/2}$ is a locally trivial (and hence trivializable) fibre bundle over entire $D$. Therefore any closed loop $l$ in $D \setminus 0$ with basepoint $\lambda_0$ defines local monodromy and variation operators

$$M^{(l)} : \tilde{H}_*(A \cap X_{\lambda_0}) \to \tilde{H}_*(A \cap X_{\lambda_0}),$$

$$\mathrm{Var}^{(l)} : \tilde{H}_*(A \cap X_{\lambda_0}, \partial B) \to \tilde{H}_*(A \cap X_{\lambda_0}), \tag{1}$$

cf. [Pham], [AVGL 1, 2], [GM 86].

The variation operator (1) is the main subject in the following *localization principle* for the global monodromy action. Suppose that $\tilde{f}$ is defined in entire $\mathbf{C}^n$ and the varieties $\Upsilon_\lambda \equiv \tilde{f}^{-1}(\lambda)$, $\lambda \in D$, are transversal to $A$ outside $B$ (and satisfy natural regularity conditions at infinity). Then the global monodromy operator $\tilde{H}_*(A \cap \Upsilon_{\lambda_0}) \to \tilde{H}_*(A \cap \Upsilon_{\lambda_0})$ defined by the loop $l$ equals the identity operator plus the composition of three operators

$$\tilde{H}_*(A \cap \Upsilon_{\lambda_0}) \to \tilde{H}_*(A \cap \Upsilon_{\lambda_0} \cap B, \partial B) \to \tilde{H}_*(A \cap \Upsilon_{\lambda_0} \cap B) \to \tilde{H}_*(A \cap \Upsilon_{\lambda_0}), \tag{2}$$

the first of which is induced by the reduction modulo $A \setminus B$, the second is the local variation operator (1), and the third is induced by the identical imbedding.

For an application, consider the problem of the integral geometry indicated in the Introduction. Let $A$ be a projective complex stratified variety and $\mathrm{reg}(A)$ the set of all hyperplanes $\Upsilon \in (\mathbf{CP}^n)^*$ transversal to the stratification of $A$. The group $\pi_1(\mathrm{reg}(A))$ is generated by finitely many simple loops embracing the surface $(\mathbf{CP}^n)^* \setminus \mathrm{reg}(A)$ close to the nonsingular points of it corresponding to the planes *simply tangent* to some strata of $A$. The condition of simple tangency of a hyperplane $\Upsilon_0$ and a stratum $\sigma$ at the point $a \in \sigma$ is equivalent to the following property: the restriction to $A$ of the (local) affine linear form distinguishing $\Upsilon_0$ has a Morse singular point at $a$. Thus the (global) monodromy action of the simple loop $l$ on the space $H_*(A \cap \Upsilon)$ (for $\Upsilon \in \mathrm{reg}(A)$ sufficiently close to $\Upsilon_0$) can again be reduced to the local variation operator (1).

In this article we show that (and how) the calculation of the operator (1) can be reduced to the calculation of similar operators for the transversal slice $A \cap T \subset T$, where $T \subset \mathbf{C}^n$ is a complex plane of complementary dimension $n - k$ transversal to $\sigma$ at $a$.

We shall need several examples of such variation operators corresponding to the 0-dimensional strata $\sigma$.

**Example 1.** Let $n = 2$ and let $A$ be the union of $s$ irreducible germs of curves of multiplicities $u_1, \ldots, u_s$ at 0. Then the variety $A \cap X_\lambda$, $\lambda \neq 0$, consists of $u \equiv u_1 + \cdots + u_s$ points, and the monodromy operator acts by cyclic permutations in any of these $s$ groups of points. Thus the operator (1), defined by the loop generating $\pi_1(D \setminus 0)$, has an $(s-1)$-dimensional kernel generated by all these groups of $u_i$ points modulo the union of all $u$ points.

**Example 2.** Let 0 be a Morse singular point of $A$, i.e., $A$ is distinguished by the condition $\Phi \equiv z_1^2 + \cdots + z_n^2 = 0$ in some local coordinates; let $\tilde{f}$ be a generic linear form. For any $\lambda \in D$ denote by $\varphi_\lambda$ the restriction of $\Phi$ on the manifold $X_\lambda \equiv \tilde{f}^{-1}(\lambda) \cap B$. Then $\varphi_0$ is a Morse function of $n - 1$ variables, and the set of functions $\varphi_\lambda$ is a deformation of it. In particular, for $\lambda \in D \setminus 0$ the variety $A \cap X_\lambda$ is smooth and $\tilde{H}_{n-2}(A \cap X_\lambda) \sim \mathbf{Z} \sim \tilde{H}_{n-2}(A \cap X_\lambda, \partial B)$. Consider the map of the base $D$ of the deformation $\{\varphi_\lambda\}$ into the base of the standard (one-dimensional) versal deformation of $\varphi_0$, inducing from it a deformation equivalent to $\{\varphi_\lambda\}$, see f.i. [AVG 1], [AVGL 1]. This map sends the basic loop $l \in \pi_1(D \setminus 0)$ into a loop which turns around the point 0 twice. Hence for even $n$ the corresponding operator $\mathrm{Var}^{(l)}$ is trivial and for odd $n$ the operator $\mathrm{Var}^{(rl)}$, defined by the $r$-fold iteration of $l$, maps a generator of the group $\tilde{H}_{n-2}(A \cap X_{\lambda_0}, \partial B)$ into $2r$ times a generator of $\tilde{H}_{n-2}(A \cap X_{\lambda_0})$.

We shall often use the following well-known facts.

**Proposition 1** (see f.i. [GM 86]). *For any point $a$ of an analytic subset $A \subset \mathbf{C}^n$ there exists a small disc $B$ centred at $a$ such that*

*a) the pair $(A \cap B, A \cap \partial B)$ is homeomorphic to the cone over $A \cap \partial B$; this homeomorphism is identical on $A \cap \partial B$ and maps $a$ into the vertex of the cone;*

*b) the intersection homology group $IH_i(A \cap B)$ is trivial if $i \geq \dim_{\mathbf{C}} A$ and is isomorphic to $IH_i(A \cap B \setminus a)$ if $i < \dim_{\mathbf{C}} A$;*

*and moreover, both assertions a) and b) also hold for all smaller concentric discs.* □

## 2 Main results

### 2.1 Induction tower

Suppose that $A, a, \sigma, B, f, D$ and $X_\lambda$ are the same as in the previous section, in particular $f$ is a Morse function on $A \cap B$ with unique critical point at $a \in \sigma$. Let us take $\delta$ for the distinguished point in $D \equiv D_\delta$ and set $X \equiv X_\delta$. Denote $A \cap B$ by $A'$. Set $m = n - k$. Consider a flag of smooth complex submanifolds

$$T \equiv T^m \subset T^{m+1} \subset \cdots \subset T^n \equiv \mathbf{C}^n, \tag{3}$$

of dimensions $m, m+1, \ldots, n$ respectively, all of which intersect $\sigma$ transversally at the point $a$. Then the intersections with the strata of $A$ define a Whitney stratification on any of varieties $A' \cap T^r \subset T^r$. If the flag (3) is generic, then for any $r = m, \ldots, n$ the restriction of $\tilde{f}$ on $A \cap T^r$ is a Morse function in $B$; we shall suppose that this condition is satisfied for all $r$. Any loop $l$ in $D \setminus 0$, considered as family of planes $X_\lambda \cap T^r$, $\lambda \in l$, defines then the variation operator

$$\operatorname{Var}_r^{(l)} : \tilde{H}_*(A' \cap X \cap T^r, \partial B) \to \tilde{H}_*(A' \cap X \cap T^r), \qquad (4)$$

in particular $\operatorname{Var}_n^{(l)} = \operatorname{Var}^{(l)}$.

Almost everywhere below we consider only the operators (4) defined by the loop $l = \{\delta \cdot e^{it}\}, t \in [0, 2\pi]$, and denote these operators simply by $\operatorname{Var}_r$.

Denote the groups participating in these operators as follows:

$$\bar{\mathcal{H}}_*(r) \equiv \tilde{H}_*(A' \cap X \cap T^r, \partial B), \quad \mathcal{H}_*(r) \equiv \tilde{H}_*(A' \cap X \cap T^r). \qquad (5)$$

Similar notation for intersection homology groups (of middle perversity) is

$$I\bar{\mathcal{H}}_*(r) \equiv I\tilde{H}_*(A' \cap X \cap T^r, \partial B), \quad I\mathcal{H}_*(r) \equiv I\tilde{H}_*(A' \cap X \cap T^r). \qquad (6)$$

Let $d$ be the complex dimension of $A \cap X \cap T^m$, so that $\dim A = d + k + 1$. Since all $A' \cap X \cap T^r$ are Stein spaces, the following statement holds for any $r = m, \ldots, n$.

**Proposition 2** (see f.i. [GM 86]). $\mathcal{H}_i(r) = I\mathcal{H}_i(r) = 0$ for $i > d + r - m$, $I\bar{\mathcal{H}}_i(r) = 0$ for $i < d + r - m$. □

## 2.2 Stabilization formulae for standard homology

**Theorem 1.** *For any $l = 1, \ldots, k$, any integer $i$ and any coefficient group there are canonical isomorphisms*

$$\mathcal{H}_{i+l}(m+l) \simeq \mathcal{H}_i(m), \qquad (7)$$

$$\bar{\mathcal{H}}_{i+l}(m+l) \simeq \mathcal{H}_i(m) \oplus \tilde{H}_{i-l+2}(A' \cap T^m, \partial B), \qquad (8)$$

*or, equivalently, for any $r = m, \ldots, n-1$ and any integer $j$*

$$\mathcal{H}_{j+1}(r+1) \simeq \mathcal{H}_j(r), \qquad (9)$$

$$\bar{\mathcal{H}}_{j+1}(r+1) \simeq \mathcal{H}_j(r) \oplus \tilde{H}_{j+1}(A' \cap T^r, \partial B); \qquad (10)$$

*the last summand in this expression (10) is equal also to the semidirect sum*

$$[\mathcal{H}_j(r)/\operatorname{Im} \operatorname{Var}_r \bar{\mathcal{H}}_j(r)] \vdash \operatorname{Ker} \operatorname{Var}_r \bar{\mathcal{H}}_{j-1}(r). \qquad (11)$$

The exact construction of these isomorphisms will be specified in § 3. More precisely, there are two different realizations of the isomorphisms (8) and (10): one of them does not depend on the coefficient group, and the other (more convenient) is admissible only if the group of coefficients allows division by 2.

For any element $\alpha \in \mathcal{H}_j(r)$ denote by $\Sigma(\alpha)$ the element of the group $\mathcal{H}_{j+1}(r+1)$ corresponding to it by the isomorphism (9), by $\rightleftharpoons(\alpha)$ the element of $\bar{\mathcal{H}}_{j+1}(r+1)$ corresponding to the element $\alpha$ of the first summand of the right-hand side of (10) in the case of coefficients divisible by 2, and by $\rightharpoonup(\alpha)$ the similar element for general coefficient group.

**Proposition 3.** *If $\alpha \in \operatorname{Im} \operatorname{Var}_r \bar{\mathcal{H}}_j(r)$ and the group of coefficients allows division by 2, then $\rightleftharpoons (\alpha) = -\rightharpoonup (\alpha)$.*

**Theorem 2.** *If the group of coefficients allows division by 2, then the obvious map*

$$J_{r+1} : \mathcal{H}_{j+1}(r+1) \to \bar{\mathcal{H}}_{j+1}(r+1) \tag{12}$$

*(reduction mod $\partial B$) acts into the first term on the right-hand side of (10) and maps any element $\Sigma(\alpha)$ into $-\rightleftharpoons (2\alpha + \operatorname{Var}_r \circ J_r(\alpha))$. In any case, $J_{r+1}$ maps $\Sigma(\alpha)$ into $-\rightharpoonup (2\alpha + \operatorname{Var}_r \circ J_r(\alpha)) - \{\text{the class of the element } \alpha \text{ in the first term of (11)} \}$.*

**Theorem 3.** *The operator $\operatorname{Var}_{r+1}$ is equal to zero on the second term of the right-hand side of (8) or (10) and maps any element $\rightleftharpoons (\alpha)$ or $\rightharpoonup (\alpha)$ of the first term into $\Sigma(\alpha)$.*

**Corollary 1.** *The first summand in (11) (and hence the additional term in the last assertion of Theorem 2) can be nontrivial only if $r = m$.* □

**Corollary 2.** *If the group of coefficients allows division by 2, then for any $l = 1, \ldots, k$, any natural $i$ and any class $\alpha \in \mathcal{H}_i(m)$ the natural homomorphism $J_{m+l} : \mathcal{H}_{i+l}(m+l) \to \bar{\mathcal{H}}_{i+l}(m+l)$ maps the element $\Sigma^l(\alpha)$ into $\rightleftharpoons \circ \Sigma^{l-1} \circ \operatorname{Var}_m \circ J_m(\alpha)$ if $l$ is even and into $-\rightleftharpoons \circ \Sigma^{l-1}(2\alpha + \operatorname{Var}_m \circ J_m(\alpha))$ if $l$ is odd.*

This follows from Theorems 2 and 3 by the induction over $l$. □

**Definition.** For any $r = m, \ldots, n-1$, the operator $\operatorname{stab}_r$ is the composite map

$$\rightharpoonup \circ \operatorname{Var}_r : \bar{\mathcal{H}}_*(r) \to \bar{\mathcal{H}}_{*+1}(r+1). \tag{13}$$

The previous theorems imply immediately the following assertions.

**Corollary 3.** $\operatorname{Var}_{r+1} \circ \operatorname{stab}_r \equiv \Sigma \circ \operatorname{Var}_r$. □

**Theorem 4.** *Suppose that for some $r \geq m$ the operator*

$$\operatorname{Var}_r : \bar{\mathcal{H}}_*(r) \to \mathcal{H}_*(r)$$

*is an isomorphism. Then for all $s \geq m$ the operators $\operatorname{Var}_s$ are also isomorphisms, and the operators $\operatorname{stab}_s$ for $s = m, \ldots, n-1$ are isomorphisms of degree 1.* □

**Corollary 4 (periodicity theorem).** *Under the assumptions of the previous theorem, if $r \leq n-2$, then two homomorphisms $\Sigma^2 : \mathcal{H}_*(r) \to \mathcal{H}_{*+2}(r+2)$ and $\operatorname{stab} \circ \operatorname{stab} : \bar{\mathcal{H}}_*(r) \to \bar{\mathcal{H}}_{*+2}(r+2)$ are isomorphisms and commute with the operators $\operatorname{Var} : \bar{\mathcal{H}}_\cdot \to \mathcal{H}_\cdot$ (which are also isomorphisms) and the obvious maps $J : \mathcal{H}_\cdot \to \bar{\mathcal{H}}_\cdot$.* □

**Corollary 5**: *the standard Picard–Lefschetz formula.*

This formula describes the variation operator (1) in the special case when $A$ is a smooth hypersurface in $\mathbf{C}^n$ and $\sigma = A$. For $n = 2$ this formula is obvious (as well as the fact that the groups $\mathcal{H}_i(n)$ and $\bar{\mathcal{H}}_i(n)$ are trivial for $i \neq n-2$ and free cyclic for $i = n-2$): indeed, in this case the monodromy action is the permutation of two points. For arbitrary $n$ this fact and the Picard–Lefschetz formula follow from this case and Theorems 1–4. □

**Example**

Under the assumptions of Theorems 1–3 let the codimension $m \equiv n - k$ of the stratum $\sigma$ be equal to 2; cf. Figure 1. The transversal slice $A' \cap T^m$ of $A$ at the point $a$ is a complex curve, and the set $A' \cap X \cap T^m$ consists of finitely many points; the number of them is the local multiplicity of this plane curve, and the monodromy and variation operators for this slice are described in Example 1 of § 1.

Let $\gamma_1$ be any of these points and $\nu$ the local multiplicity at $a$ of the local irreducible component of $A' \cap T^m$ containing $\gamma_1$. Denote by $\gamma_2, \ldots, \gamma_\nu$ the points obtained consecutively from $\gamma_1$ by the iterations of the monodromy of the plane $X_\xi$ along the loop $l = \{\xi = \delta \cdot e^{i\tau}, \tau \in [0, 2\pi]\}$. Obviously, all these $\nu$ points are different, and the monodromy operator acts on them by cyclic permutation.

For any $j = 1, \ldots, \nu$, denote by $\Gamma_j$ the element $\operatorname{stab}^{n-2}(\gamma_j) \in \bar{\mathcal{H}}_{n-2}(n)$.

**Proposition 4.** *Suppose that the curve $A' \cap T^2$ is locally irreducible at the point $a$, so that its multiplicity is equal to $\nu$. Then the group $\bar{\mathcal{H}}_{n-2}(n)$ is isomorphic to $\mathbf{Z}^{\nu-1}$ and is generated by the elements $\Gamma_j$, $j = 1, \ldots, \nu$, subject to the relation $\Gamma_1 + \cdots + \Gamma_\nu = 0$. The group $\mathcal{H}_{n-2}(n)$ is also isomorphic to $\mathbf{Z}^{\nu-1}$ and is generated by the elements $\Sigma^{(n-2)}(\gamma_{j+1} - \gamma_j)$ subject to the obvious relation.*

*If $n$ is even, then $\operatorname{Var}_n^{(\nu \cdot l)} \equiv 0$, i.e. the $\nu$-fold transportation of any element $\Gamma_j$ along the loop $l$ maps this element into itself and adds nothing to the group $\mathcal{H}_{n-2}(n)$.*

*If $n$ and $\nu$ are odd, then the same is true for the $2\nu$-fold iteration of the monodromy: $\operatorname{Var}_n^{(2\nu \cdot l)} \equiv 0$.*

*If $n$ is odd and $\nu$ is even, then for any $j = 1, \ldots, \nu$*

$$\operatorname{Var}_n^{(\nu \cdot l)}(\Gamma_j) = \Sigma^{n-2}((\gamma_{j+1} - \gamma_j) + (\gamma_{j+3} - \gamma_{j+2}) + \ldots + (\gamma_{j-1} - \gamma_{j-2})) \tag{14}$$

*(where the lower indices in the right-hand expression are considered as residues mod $\nu$).*

*The cycle (14) defines a zero element of the relative homology group $\bar{\mathcal{H}}_{n-2}(n)$.*

This is a direct corollary of Theorems 1–3. □

**Remark.** The formula (14) is essential in the integral geometry: e.g., it allows to prove the logarithmic ramification of volume functions close to the planes tangent to the cuspidal edges in the odd-dimensional Newton's problem on algebraically integrable bodies, see [V].

## 2.3 Stabilization of the intersection homology

We consider only the intersection homology groups (of middle perversity) with the coefficients in a field of characteristic zero. (In fact, we use only the assumptions that all groups $I\mathcal{H}_i(m)$, $I\bar{\mathcal{H}}_i(m)$ and $I\mathcal{H}_d(m)/\operatorname{Im}\operatorname{Var}_m(I\bar{\mathcal{H}}_d(m))$ are torsion-free and for any $i = 1, \ldots, d$ the groups $I\mathcal{H}_i(m)$ and $I\bar{\mathcal{H}}_{2d-i}(m)$ are dual over the group of coefficients with respect to the intersection index; here $d$ and $m$ are the same as in § 2.1.)

**Theorem 5.** *For any $l = 1, \ldots, k$, there are natural isomorphisms*

$$I\mathcal{H}_{d+l}(m+l) \simeq \operatorname{Im} \operatorname{Var}_m(I\bar{\mathcal{H}}_d(m)), \tag{15}$$

$$I\bar{\mathcal{H}}_{d+l}(m+l) \simeq I\bar{\mathcal{H}}_d(m)/\text{Ker Var}_m, \tag{16}$$

$$I\mathcal{H}_i(m+l) \simeq I\mathcal{H}_i(m) \text{ for } i \leq d-1, \tag{17}$$

$$I\bar{\mathcal{H}}_j(m+l) \simeq I\bar{\mathcal{H}}_{j-2l}(m) \text{ for } j \geq d+2l+1, \tag{18}$$

$$I\mathcal{H}_d(m+l) \simeq I\mathcal{H}_d(m)/\text{Im Var}_m(I\bar{\mathcal{H}}_d(m)), \tag{19}$$

$$I\bar{\mathcal{H}}_{d+2l}(m+l) \simeq \text{Ker Var}_m(I\bar{\mathcal{H}}_d(m)), \tag{20}$$

$$I\mathcal{H}_i(m+l) \simeq 0 \text{ for } i \in [d+1, d+l-1] \text{ or } i > d+l, \tag{21}$$

$$I\bar{\mathcal{H}}_j(m+l) \simeq 0 \text{ for } j \in [d+l+1, d+2l-1] \text{ or } i < d+l. \tag{22}$$

In particular, all groups $I\mathcal{H}_{d+l}(m+l)$ (respectively, $I\bar{\mathcal{H}}_{d+l}(m+l)$) with $l \geq 1$ are canonically isomorphic to one another.

**Remark.** In general, the requirement $l \geq 1$ in the last statement cannot be extended to $l = 0$: indeed, it follows from Example 2 of § 1 that the operator $\text{Var}_m : I\bar{\mathcal{H}}_d(m)) \to I\mathcal{H}_d(m))$ can be not isomorphic.

**Remark.** The formulae (17) and (18) also could be rewritten in the uniform way ((19), (20)) with the dimension $d$ replaced by $i$ and $j - 2l$ respectively: indeed, in these dimensions the operator $\text{Var}_m$ is obviously trivial.

For arbitrary element $\alpha$ of the right-hand group in formula (17) or (19) (respectively, (18) or (20)) denote by $I_*^l(\alpha)$ (respectively, $(I^*)^{-l}(\alpha)$) the element in the left side corresponding to $\alpha$ via this formula.

For arbitrary element $\alpha$ of the group $I\bar{\mathcal{H}}_d(m)$ denote by $\infty^l(\alpha)$ the element of the group $I\mathcal{H}_{d+l}(m+l)$ corresponding to the class $\text{Var}_m(\alpha)$ via the isomorphism (15) and by $\text{stab}^l(\alpha)$ the element of $I\bar{\mathcal{H}}_{d+l}(m+l)$ corresponding to the coset $\{\alpha\}$ via the isomorphism (16).

**Theorem 6.** *Let $l$ be a natural number, $l \leq k$. Then for any $\alpha, \beta \in I\bar{\mathcal{H}}_d(m)$,*

$$\langle \text{stab}^l(\alpha), \infty^l(\beta) \rangle = (-1)^{l(d+(l-1)/2)} \langle \alpha, \text{Var}_m(\beta) \rangle; \tag{23}$$

*for any elements $\alpha \in I\bar{\mathcal{H}}_{2d-i}(m)$ and $\beta \in I\mathcal{H}_i(m)$, $i < d$,*

$$\langle (I^*)^{-l}(\alpha), I_*^l(\beta) \rangle = \langle \alpha, \beta \rangle, \tag{24}$$

*as well as for any $\alpha \in \text{Ker Var}_m(I\bar{\mathcal{H}}_d(m))$ and $\beta \in I\mathcal{H}_d(m)$.*

**Corollary.** *The groups $\text{Im Var}_m \subset I\mathcal{H}_d(m)$ and $\text{Ker Var}_m \subset I\bar{\mathcal{H}}_d(m)$ are orthogonal with respect to the intersection form.* $\square$

**Theorem 7.** *For any $l = 1, \ldots, k$, the homomorphism $\text{Var}_{m+l}$ maps the element $\text{stab}^l(\alpha)$ into $\infty^l(\alpha)$.*

**Theorem 8.** *For any $l = 1, \ldots, k$, the obvious homomorphism*

$$J_{m+l} : I\mathcal{H}_{d+l}(m+l) \to I\bar{\mathcal{H}}_{d+l}(m+l) \tag{25}$$

*maps any element $\infty^l(\alpha)$ into $-\text{stab}^l(2\alpha + J_m \circ \text{Var}_m(\alpha))$ if $l$ is odd and into $\text{stab}^l(J_m \circ \text{Var}_m(\alpha))$ if $l$ is even.*

Denote by $\rho_r$ and $\bar{\rho}_r$ the obvious maps of the intersection homology groups $I\mathcal{H}_{d+r-m}(r)$, $I\bar{\mathcal{H}}_{d+r-m}(r)$ into the standard homology groups $\mathcal{H}_{d+r-m}(r)$, $\bar{\mathcal{H}}_{d+r-m}(r)$.

**Theorem 9.** *For any $r = m, \ldots, n$ and any $l = 1, \ldots, n-r$, $\text{Var}_r \circ \bar{\rho}_r \equiv \rho_r \circ \text{Var}_r$; $\bar{\rho}_{r+l} \circ \text{stab}^l \equiv \text{stab}^l \circ \bar{\rho}_r$; $\rho_{r+l} \circ \infty^l \equiv \Sigma^l \circ \text{Var}_r \circ \bar{\rho}_r$.*

# 3 Stabilization of the vanishing cycles

## 3.1 Adapted coordinates and polydisc $B'$

Let $A, a, \sigma, B, \tilde{f}, f, D$ and $X_\lambda$ be the same as in the previous sections.

**Definition.** A local analytic coordinate system $\{z_1, \ldots, z_n\}$ in $\mathbf{C}^n$ with origin at $a$ is called *adapted* if $z_n \equiv \tilde{f}$, the tangent space to $\sigma$ at $a$ is spanned by the vectors $\partial/\partial z_1, \ldots, \partial/\partial z_k$ (so that the restrictions of the functions $z_1, \ldots, z_k$ constitute a local coordinate system on $\sigma$), and in restriction to $\sigma$

$$z_n \equiv z_1^2 + \cdots + z_k^2.$$

In Figure 1b a real version of this situation is shown, where $n = 3$, the plane $X \equiv X_\lambda$ is given by the equation $z_3 \equiv \lambda$, and $z_1$ is the coordinate along the stratum $\sigma$. The transversal slice of this picture by the plane $\{z_1 = 0\}$ is shown in Figure 1a.

By the Morse lemma, adapted coordinates always exist; let us fix such a coordinate system. Without loss of generality we can define the flag (3) by the conditions $T^r = \{z \mid z_1 = \cdots = z_{n-r} = 0\}$.

We can assume that $B \equiv B_\varepsilon$ is a disc of radius $\varepsilon$ with respect to the standard Hermitian metric defined by these coordinates $z_1, \ldots, z_n$. Moreover, in our considerations we can replace the disc $B$ by a polydisc defined by these coordinates. Namely, let $B' \subset B$ be the polydisc $\{z \mid |z_i| \leq \varepsilon/n \text{ for all } i\}$ and suppose that the number $\delta$ (participating in the definition of the disc $D \equiv D_\delta \subset \mathbf{C}^1$) is sufficiently small with respect to $\varepsilon$ and $\varepsilon^2$.

**Proposition 5.** *For every surface $X_\lambda \subset B$ defined by the equation $z_n \equiv \lambda$, $\lambda \in D$, the pair $(A \cap X_\lambda \cap B', A \cap X_\lambda \cap \partial B')$ is homeomorphic to the pair $(A \cap X_\lambda, A \cap X_\lambda \cap \partial B)$, and the corresponding homotopy equivalence of quotient spaces $(A \cap X_\lambda)/(A \cap X_\lambda \cap \partial B) \to (A \cap X_\lambda \cap B')/(A \cap X_\lambda \cap \partial B')$ is realized by the obvious factorization map contracting the whole complement of $B'$ to one point. Moreover, for every $r = m, \ldots, n$ and every $\lambda \in D$, the pair $(A \cap X_\lambda \cap T^r \cap B', A \cap X_\lambda \cap T^r \cap \partial B')$ is also naturally homeomorphic to $(A \cap X_\lambda \cap T^r, A \cap X_\lambda \cap T^r \cap \partial B)$.*

*Proof.* It is sufficient to consider only the case $r = n$. First we prove our assertion for the plane $X_0$. Consider two metrics in $\mathbf{C}^n$, $|z| = (\sum |z_i|^2)^{1/2}$ and $\|z\| = \max |z_i|$.

**Lemma 1.** *If the radius $\varepsilon$ of the disc $B$ is sufficiently small, then there exists a homeomorphism of the analytic variety $A \cap X_0 \subset B$ onto the cone over its boundary $A \cap X_0 \cap \partial B$ which is identical on this boundary and maps the point $a$ into the vertex of the cone, while both functions $|z|, \|z\|$ strictly increase along the images of all generating segments of the cone under this homeomorphism.*

This is essentially a special case of Proposition 1a: the only additional requirement about the increase of $|z|$ and $\|z\|$ is ensured by the proof of it, see f.i. [Milnor]. □

The assertion of Proposition 5 for the plane $X_0$ follows immediately from this lemma.

Let $\mathcal{B}$ be an open disc centred at $a$ and contained in $B'$. Since the plane $X_0$ is transversal to the stratified variety $A$ everywhere in $B \setminus a$, then any plane $X_\lambda$ with very small $|\lambda|$ is also transversal to $A$ in $B \setminus \mathcal{B}$. Hence, by Thom's isotopy theorem, all the varieties $A \cap X_\lambda \cap (B \setminus \mathcal{B})$ with sufficiently small $|\lambda|$ are isotopic to $A \cap X_0 \cap (B \setminus \mathcal{B})$; moreover, for every generating segment of the cone $A \cap X_0 \simeq C(A \cap X_0 \cap \partial B)$ the part of this segment belonging to $B \setminus \mathcal{B}$ goes under this isotopy into a segment along which both functions $\|z\|$ and $|z|$ are again monotonic. This implies Proposition 5. □

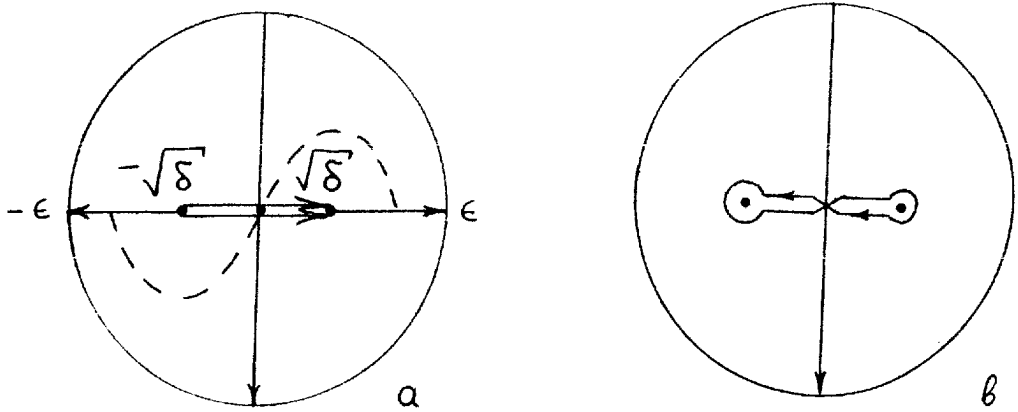

Figure 2: The base $\Omega$ of the fibre bundle $z_{n-r}$

Thus, everywhere in the proof of Theorems 1–8 we can replace the disc $B$ by the polydisc $B'$, in particular $A'$ and $X$ are the notation for $A \cap B'$ and $\tilde{f}^{-1}(\delta) \cap B'$ respectively.

## 3.2 Stabilization of the standard homology groups

For any $r = m, \ldots, n$ denote the $2r$-dimensional polydisc $B' \cap T^r$ by $B^{(r)}$ and the characteristic radius $\varepsilon/n$ of all such polydiscs by $\epsilon$. Remember the notation $X \equiv X_\delta$. Denote the $\epsilon$-disc $\{z \in \mathbf{C}^1 \mid |z| \le \epsilon\}$ by $\Omega$. For arbitrary $r < n$ consider the projection

$$z_{n-r} : A \cap X \cap B^{(r+1)} \to \Omega. \qquad (26)$$

For any $t \in \Omega$ denote by $\mathcal{F}_t$ the fibre $A \cap X \cap B^{(r+1)} \cap \{z \mid z_{n-r} = t\}$ of this projection, and by $\partial \mathcal{F}_t$ the boundary $\mathcal{F}_t \cap \partial B^{(r+1)}$ of this fibre. In particular, $\mathcal{F}_0 \equiv A \cap X \cap B^{(r)}, \tilde{H}_*(\mathcal{F}_0, \partial \mathcal{F}_0) \equiv \bar{\mathcal{H}}_*(r)$ and $\tilde{H}_*(\mathcal{F}_0) \equiv \mathcal{H}_*(r)$.

Each of the fibres $\mathcal{F}_{\sqrt{\delta}}$, $\mathcal{F}_{-\sqrt{\delta}}$ contains a distinguished point, the critical point of the restriction of $z_{n-r}$ to the manifold $\sigma \cap X \cap T^{r+1}$. (These points can however be nonsingular points of the whole $A' \cap X \cap T^{r+1}$.)

**Lemma 2.** *If the radius $\epsilon$ of the polydisc $B'$ is sufficiently small and $\delta$ is sufficiently small with respect to $\epsilon^2$, then the projection (26) is a locally trivial fibre bundle over the disc $\Omega$ with two points $\pm\sqrt{\delta}$ removed, the standard fibre of which is homeomorphic to the variety $\mathcal{F}_0$, and the fibres over the exceptional points $\sqrt{\delta}$ and $-\sqrt{\delta}$ are homeomorphic to* the cones over $\partial\mathcal{F}_0$ with the distinguished points as vertices; in particular these fibres $\mathcal{F}_{\pm\sqrt{\delta}}$ are homotopy trivial; see Figures 2, 1.

This lemma follows directly from the construction, from Thom's isotopy theorem and from Proposition 1a. □

Denote by $W$ the union of two very small (with respect to $\epsilon$ and $\delta$) discs in $\mathbf{C}^n$ around the distinguished points and by $\Im$ the intersection of $A' \cap X \cap T^{r+1}$ with the pre-image of the imaginary axis $\{z \in \Omega \mid \operatorname{Re} z = 0\}$ under the projection (26).

**Proposition 6.** *The quotient space $\Im/(\Im \cap \partial B')$ is a deformation retract of the quotient space $(A' \cap X \cap T^{r+1} \setminus W)/(A' \cap X \cap T^{r+1} \cap \partial B')$. In particular, for any $i$ the group $\tilde{H}_{i+1}((A' \cap X \cap T^{r+1} \setminus W), \partial B')$ is naturally isomorphic to $\bar{\mathcal{H}}_i(r)$.*

Indeed, the centrifugal vector fields issuing from these two distinguished points push the points of $A' \cap X \cap T^{r+1} \setminus W$ out from the small neighborhoods of the fibres $\mathcal{F}_{\sqrt{\delta}}$, $\mathcal{F}_{-\sqrt{\delta}}$, so that their images after this deformation intersect these neighborhoods only in points of $\partial B^{(r+1)}$, which can be ignored. Then, using the fibre bundle structure over the complement of the points $\pm\sqrt{\delta}$ in $\Omega$, we push this complement into the pre-image of the imaginary axis. The last (homological) assertion of the proposition follows from the fact that over this axis the bundle (26) is trivializable. □

Consider two variation operators $V_{+,-} : \tilde{H}_*(\mathcal{F}_0, \partial\mathcal{F}_0) \to \tilde{H}_*(\mathcal{F}_0)$ defined by the simple loops in $\Omega$ corresponding to the segments $[0, \sqrt{\delta}]$ and $[0, -\sqrt{\delta}]$. (For instance the $\infty$-shaped loop in Figure 2b is a composition of the second of these loops and the loop inverse to the first of them.)

**Lemma 3.** *The operators $V_+$, $V_-$ are equal to each other and to the operator $\operatorname{Var}_r : \bar{\mathcal{H}}_*(r) \to \mathcal{H}_*(r)$. The same is true also for three similar operators $I\bar{\mathcal{H}}_*(r) \to I\mathcal{H}_*(r)$.*

Indeed, all three operators are defined by three homotopic loops in the space of $(r-1)$-dimensional complex planes parallel to $X \cap T^r$ and transversal to $A$. □

Let $\alpha$ be any element of the group $\mathcal{H}_j(r)$. Using the fibre bundle structure (26), we transport the realizing $\alpha$ cycle in $\mathcal{F}_0$ over the S-shaped path in Figure 2a to some interior points of the segments $[\sqrt{\delta}, \epsilon]$ and $[-\epsilon, -\sqrt{\delta}]$. Then, transporting the resulting cycles over these segments, we sweep out two $(j+1)$-dimensional chains in $A' \cap X \cap B^{(r+1)}$. We orient them by the pair of orientations, the first of which is the orientation of the base segment chosen as shown in Figure 2a, and the second is induced by the original orientation of our cycle $\alpha$ over the basepoint 0. Since the fibres over the points $\pm\sqrt{\delta}$ are homotopically trivial, we can add to these chains two chains in these singular fibres, which span their boundaries, so that the resulting two chains become relative cycles (modulo $\partial B^{(r+1)}$).

If the basic coefficient group admits division by 2, then the class $\rightleftharpoons (\alpha) \in \bar{\mathcal{H}}_{j+1}(r+1)$ corresponding to $\alpha$ in the first term of the right-hand side of (10) is defined by half the difference of these two cycles; $\rightharpoonup(\alpha)$ is just the class of the first of them.

The sum of these two cycles is of course a cycle too. This cycle is homologous to zero if and only if $\alpha \in \operatorname{Im}\operatorname{Var}_r \bar{\mathcal{H}}_j(r)$, and we get the first term of (11). If the relative cycle $\gamma \in \bar{\mathcal{H}}_{j-1}(r)$ belongs to the subgroup $\operatorname{Ker}\operatorname{Var}_r$, then it can be transported continuously into all the fibres of our fibration over $\Omega$, including the degenerate fibres (cf. the proof of Theorem 12 below); this two-dimensional family of $(j-1)$-dimensional cycles forms an $(j+1)$-dimensional cycle, which is included in the second term of (11).

The class $\Sigma(\alpha) \in \mathcal{H}_{j+1}(r+1)$ realizing the isomorphism (9) is obtained from $\alpha$ by the suspension operation: it is swept out by an appropriate family of cycles obtained from $\alpha$ by transport over the segment $[-\sqrt{\delta}, \sqrt{\delta}]$ and subsequent contraction of the boundary of the resulting cycle inside the fibres over the endpoints $-\sqrt{\delta}$ and $\sqrt{\delta}$.

The operation $\Sigma^l$ realizing the formula (7) is just the $l$-th iteration of $\Sigma$; the operations $\rightleftharpoons^l$ or $\multimap^l$, realizing the first part of (8), are defined as $\rightleftharpoons \circ \Sigma^{l-1}$ or $\multimap \circ \Sigma^{l-1}$ respectively.

**Theorem 10.** *The operator* stab *(see (13)) can be realized by the relative cycles* $\mathrm{stab}(\alpha)$ *swept out by the corresponding cycles* $\alpha \in \bar{\mathcal{H}}_*(r)$ *in transport over the imaginary axis oriented downwards; see Figure 2.*

## 3.3 Stabilization of the intersection homology

The operators $I_*^l$ realizing the formulae (17), (19) are induced by the identical imbeddings $A' \cap X \cap T^m \hookrightarrow A' \cap X \cap T^{m+l}$. The operator $(I^*)^l : I\bar{\mathcal{H}}_{2d+2l-i} \to I\bar{\mathcal{H}}_{2d-i}$ is dual to $I_*^l$; the formulae (18), (20) assert that this dual operator is isomorphic for $i < d$, and for $i = d$ it is monomorphic and its image is equal to $\mathrm{Ker}\,\mathrm{Var}_m$.

Formula (24) follows immediately from this realization.

For any $r = m, \ldots, n-1$ and $j$ the operator stab $: I\bar{\mathcal{H}}_j(r) \to I\bar{\mathcal{H}}_{j+1}(r+1)$ is defined as similar operator in Theorem 10: for any element $\alpha \in I\bar{\mathcal{H}}_j(r)$ a cycle representing the class $\mathrm{stab}(\alpha)$ is swept out by any cycle representing $\alpha$ transported over the imaginary axis.

The operator $\mathrm{stab}^l$ realizing formula (16) is just the $l$-th power of stab.

**Theorem 11.** *For any natural* $l \leq k$, *the map*

$$stab : I\bar{\mathcal{H}}_{d+l-1}(m+l-1) \to I\bar{\mathcal{H}}_{d+l}(m+l) \tag{27}$$

*is epimorphic.*

*Proof.* For any finite subset of a complex variety of dimension $q \geq 1$, any intersection homology class of dimension $\leq q$ of this variety can be realized by a cycle avoiding this subset. In particular, any element of the right-hand group in (27) can be realized as a relative cycle in $A' \cap X \cap T^{m+l} \setminus W$ (mod $\partial B'$). Then we deform this cycle into the space $\Im \cup \partial B'$ as in the proof of Proposition 6 and use the Künneth formula for the relative intersection homology of $\Im$, see [GM 86], § II.6.1. $\square$

Finally, the operator $\infty : I\bar{\mathcal{H}}_j(r) \to I\bar{\mathcal{H}}_{j+1}(r+1)$ is defined as follows. Consider again the projection (26) with fibres $\mathcal{F}_t$. For any $\epsilon' < \epsilon$ sufficiently close to $\epsilon$ the restriction of this projection on the complement of the polydisc $B''$ of radius $\epsilon'$ in $\mathbf{C}^n$, $z_{n-r} : A \cap X \cap T^{r+1} \cap (B'_\epsilon \setminus B''_{\epsilon'}) \to \Omega_{\epsilon'}$, is trivializable over $\Omega_{\epsilon'}$, in particular there is a homeomorphism

$$\Lambda : A \cap X \cap T^{r+1} \cap (B'_\epsilon \setminus B''_{\epsilon'}) \cap z_{n-r}^{-1}(\Omega_{\epsilon'}) \sim [\mathcal{F}_0 \cap (B'_\epsilon \setminus B''_{\epsilon'})] \times \Omega_{\epsilon'}. \tag{28}$$

Denote by $\Lambda_1$ and $\Lambda_2$ compositions of $\Lambda$ with projections onto the factors in the right-hand side of (28); since (28) is a trivialization of the fibre bundle (26), $\Lambda_2 \equiv z_{n-r}$.

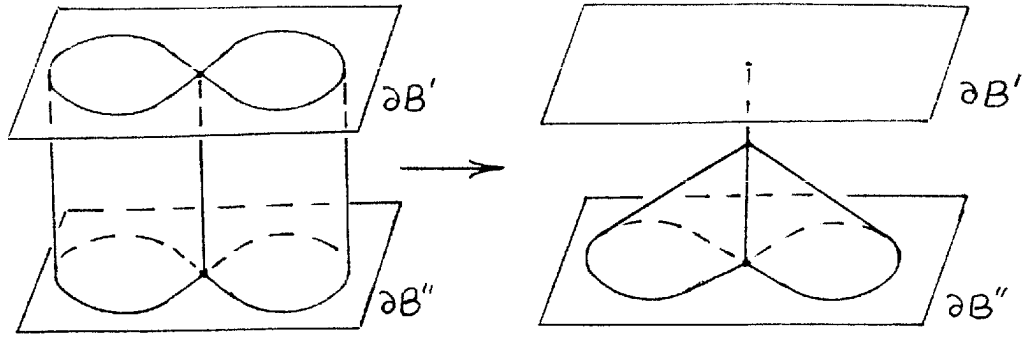

Figure 3: Regularization of the cycle $\tilde{\infty}(\alpha)$

We can suppose that the $\infty$-shaped path in $\Omega$ lies inside $\Omega_{\epsilon'}$, see Figure 2b. For any class $\gamma \in I\bar{\mathcal{H}}_j(r)$ we transport a relative cycle realizing $\gamma$ over this path in such a way that outside $B''_{\epsilon'}$ this transportation is compatible with the trivialization (28). Let $\tilde{\infty}(\gamma)$ be the cycle swept out by this transportation. The desired class $\infty(\gamma)$ is the class of a "regularization" of $\tilde{\infty}(\gamma)$, i.e., of an absolute cycle having the same class in the group $I\bar{\mathcal{H}}_{j+1}(r+1)$.

Indeed, let $h : \Omega \times [0,1] \to \Omega$ be a homotopy contracting the disc $\Omega$ to the point $0$, $h(\cdot, 0) \equiv \mathrm{id}$, $h(\Omega, 1) = \{0\}$. Using the trivialization (28) we lift this homotopy to the cycle $\tilde{\infty}(\gamma) \cap (B'_\epsilon \setminus B''_{\epsilon'})$ in such a way that $\tilde{\infty}(\gamma)$ becomes a compact cycle, see Figure 3. Namely, we define the map $\chi : \tilde{\infty}(\gamma) \to A \cap X \cap B^{(r+1)}$ by the following condition. For any $z \in \tilde{\infty}(\gamma)$, $\chi(z)$ is equal to

- $z$  if $\|z\| \leq \epsilon'$,
- $\Lambda^{-1}(\Lambda_1(\Lambda(z)) \times \{0\})$  if $\|z\| \geq (\epsilon + \epsilon')/2$,
- $\Lambda^{-1}(\Lambda_1(\Lambda(z)) \times h(\Lambda_2(\Lambda(z)), ((\|z\| - \epsilon') \cdot 2/(\epsilon - \epsilon'))))$  if $\epsilon' \leq \|z\| \leq (\epsilon + \epsilon')/2$.

The "collar" part of the cycle $\chi(\tilde{\infty}(\gamma))$ that lies in the set where $\|z\| > (\epsilon + \epsilon')/2$ is obviously a zero chain. Removing this part we obtain a compact cycle; the desired element $\infty(\gamma)$ is defined as the class of this cycle.

Define the operator $\infty^l$ realizing the isomorphism (15) as the composition $\infty \circ \mathrm{stab}^{l-1}$.

**Lemma 4.** *In a neighbourhood of the set $\Im$ (see Proposition 6) a cycle realizing $\infty(\gamma)$ coincides with the cycle $\mathrm{Var}_r(\gamma) \subset \mathcal{F}_0$ transported over the real segment $(-\sqrt{\delta}, \sqrt{\delta})$ via the fibre bundle (26).*

This follows immediately from the construction and from Lemma 3.  □

**Corollary 1.** *For any $l \geq 1$ and $\alpha, \beta \in I\bar{\mathcal{H}}_{d+l-1}(m+l-1)$,*

$$\langle \mathrm{stab}(\alpha), \infty(\beta) \rangle = (-1)^{d+l-1} \langle \alpha, \mathrm{Var}_{m+l-1}(\beta) \rangle. \quad \square$$

**Corollary 2.** $\dim I\bar{\mathcal{H}}_{d+l}(m+l) \geq \dim \mathrm{Im}\, \mathrm{Var}_{m+l-1}(I\bar{\mathcal{H}}_{d+l-1}(m+l-1))$.

This follows immediately from the previous corollary and the fact that the intersection pairing $I\bar{\mathcal{H}}_{d+l-1}(m+l-1) \otimes I\bar{\mathcal{H}}_{d+l-1}(m+l-1) \to \mathbf{C}$ is nondegenerate. □

**Theorem 12.** *For any $r = m, \ldots, n-1$ and any $j$, if $\gamma \in \mathrm{Ker\,Var}_r \subset I\bar{\mathcal{H}}_j(r)$, then $\mathrm{stab}(\gamma) = 0$. In particular, $\mathrm{stab}(\gamma) = 0$ if $j \neq d + r - m \equiv \dim_{\mathbf{C}} \mathcal{F}_0$.*

**Corollary 1.** *For any $l = 1, \ldots, k$, the dimension of the group $I\bar{\mathcal{H}}_{d+l}(m+l)$ is not greater than that of $I\bar{\mathcal{H}}_{d+l-1}(m+l-1)/\mathrm{Ker\,Var}_{m+l-1}$.*

This follows immediately from Theorems 12 and 11. □

**Corollary 2.** *For any $l = 1, \ldots, k$, $I\bar{\mathcal{H}}_{d+l}(m+l) \simeq \mathrm{Im\,Var}_{m+l-1}(I\bar{\mathcal{H}}_{d+l-1}(m+l-1))$, $I\bar{\mathcal{H}}_{d+l}(m+l) \simeq I\bar{\mathcal{H}}_{d+l-1}(m+l-1)/\mathrm{Ker\,Var}_{m+l-1}$, and these isomorphisms are realized by operations $\infty$ and $\mathrm{stab}$.*

Indeed, $\dim I\bar{\mathcal{H}}_{d+l}(m+l) \geq \dim \mathrm{Im\,Var}_{m+l-1}(I\bar{\mathcal{H}}_{d+l-1}(m+l-1)) = \dim(I\bar{\mathcal{H}}_{d+l-1}(m+l-1)/\mathrm{Ker\,Var}_{m+l-1}) \geq \dim I\bar{\mathcal{H}}_{d+l}(m+l) = \dim I\bar{\mathcal{H}}_{d+l}(m+l)$; here both equalities follow from the Poincaré duality, and the inequalities from two previous corollaries. □

*Proof of Theorem 12.* For $j < d + r - m$ this theorem is obvious (since $I\bar{\mathcal{H}}_j(r) = 0$ for such $j$); let us prove it for $j \geq d+r-m$. Let $\breve{\Omega}_-$ be the part of the halfdisc $\Omega \cap \{z \mid \mathrm{Re}\, z \leq 0\}$ not contained inside the $\infty$-shaped loop in Figure 2b. Denote by $\nu$ the radius of the small circles at the ends of this loop. Consider a two-parametric family of $j$-dimensional cycles, obtained from $\gamma$ by transportation into the points of $\breve{\Omega}_-$ via fibre bundle structure (26) and coinciding over the segment $[0, \epsilon \cdot i)$ with the family sweeping the cycle $\mathrm{stab}(\gamma)$. Arguing as in the construction of the cycle $\infty(\gamma)$ we obtain that the boundary of the resulting $(j+2)$-dimensional chain is equal to the sum of a) $\mathrm{stab}(\gamma)$, b) a one-parametric continuous family of $j$-dimensional absolute cycles in the fibres $\mathcal{F}_t$ where $t$ runs over the segments $[-\sqrt{\delta}+\nu, 0]$ and $[0, -\epsilon \cdot i)$, and c) a compact chain in the union of fibres $\mathcal{F}_\cdot$ over the disc of radius $\nu$ centred at $-\sqrt{\delta}$. Since $\gamma \in \mathrm{Ker\,Var}_r$, all $j$-dimensional cycles sweeping the part b) are homologous to zero in the corresponding fibres $\mathcal{F}_t$. After the fibre-wise continuous contraction of them, we get that $\mathrm{stab}(\gamma)$ is homologous to an absolute cycle whose projection onto $\Omega$ belongs to an arbitrarily small neighbourhood of the point $-\sqrt{\delta}$. Using Proposition 1a and the fact that the local triviality of the fibre bundle (26) over $\Omega_-$ fails only in an arbitrarily small neighbourhood of the distinguished point in the fibre $\mathcal{F}_{-\sqrt{\delta}}$, we get that this absolute cycle can be contracted into such a neighbourhood of this distinguished point, and hence, by Proposition 1b, is homologous to zero. □

# 4 Proofs of the stabilization theorems for ordinary homology

## 4.1 An exact sequence

**Lemma 5** (cf. [GM 86], § II.6.3). *Suppose that $M$ is a complex analytic subvariety in $\mathbf{C}^q$ with fixed analytic Whitney stratification, $M \ni 0$, $\Delta \equiv \Delta_\phi$ is an open disc in $\mathbf{C}^q$ of a very small radius $\phi$ centred at 0, and the restriction on $M$ of the linear function $z_1$ has an isolated singularity at 0, so that all planes $\Upsilon_\lambda \equiv \{z \mid z_1 = \lambda\}$ with sufficiently small $|\lambda|$ are* *transversal to $M$ outside $0$ in $\Delta$. Let $\psi$ be a positive number sufficiently small with respect to $\phi$; denote $\Delta \setminus 0$ by $\breve{\Delta}$ and $\Upsilon_\lambda \cap \Delta$ by $\Upsilon'_\lambda$. Then the homology groups of $M \cap \breve{\Delta}$ participate in the following exact sequences*

$$\to \tilde{H}_i(M \cap \Upsilon'_\psi) \to \tilde{H}_i(M \cap \breve{\Delta}) \to \tilde{H}_{i-1}(M \cap \bar{\Upsilon}'_\psi, \partial \bar{\Upsilon}'_\psi) \to \tilde{H}_{i-1}(M \cap \Upsilon'_\psi) \to, \qquad (29)$$

$$\to I\tilde{H}_i(M \cap \Upsilon'_\psi) \to I\tilde{H}_i(M \cap \breve{\Delta}) \to I\tilde{H}_{i-1}(M \cap \bar{\Upsilon}'_\psi, \partial \bar{\Upsilon}'_\psi) \to I\tilde{H}_{i-1}(M \cap \Upsilon'_\psi) \to, \quad (30)$$

*the first and the fourth arrows in any of which are the variation operators defined by the family of spaces $M \cap \Upsilon'_\lambda$, $\lambda = \psi e^{i\tau}, \tau \in [0, 2\pi]$.*

*Proof.* Similarly to § 3.1, we can replace $\breve{\Delta}$ by the difference of two cylinders, $\Gamma' \setminus \Gamma''$, where $\Gamma' = \{z \mid |z_1| < \psi, \|z_2, \ldots, z_q\| < \phi\}$, $\Gamma'' = \{z \mid |z_1| \le \psi - \psi^2, \|z_2, \ldots, z_q\| \le \phi - \phi^2\}$. This difference can be covered by three open subsets $\Phi_+, \Phi_-$ and $\nabla$, distinguished by conditions $\{\operatorname{Re} z_1 > -\psi^2, |z_1| > \psi - \psi^2\}$, $\{\operatorname{Re} z_1 < \psi^2, |z_1| > \psi - \psi^2\}$, and $\|z_2, \ldots, z_q\| > \phi - \phi^2$ respectively. The sets $M \cap \Phi_\pm$ (respectively, $M \cap \nabla$; respectively, the pair $(M \cap \Phi_-, M \cap \Phi_- \cap \nabla)$) are obviously homeomorphic to the direct product of an open 2-dimensional disc and $M \cap \Upsilon'_\psi$ (respectively, a collar of the boundary of $M \cap \bar{\Upsilon}'_\psi$; respectively, the pair $(M \cap \bar{\Upsilon}'_\psi,$ such a collar).

The exact sequence of the triple $(M \cap (\Gamma' \setminus \Gamma''), M \cap (\Phi_+ \cap \Phi_-), M \cap \Phi_+)$ proves that the group $\tilde{H}_{i+1}(M \cap (\Gamma' \setminus \Gamma''), M \cap \Phi_+)$ is isomorphic to $\tilde{H}_i(M \cap \bar{\Upsilon}'_\psi, \partial \bar{\Upsilon}'_\psi)$ and the boundary homomorphism of this group into $\tilde{H}_i(M \cap \Phi_+)$ coincides via this isomorphism with the operator Var. The exact sequence of the pair $(M \cap (\Gamma' \setminus \Gamma''), M \cap \Phi_+)$ is thus isomorphic to (29). The sequence (30) is justified in exactly the same way. □

**Corollary 1.** *There are exact sequences similar to (29) (respectively, (30)), in which instead of the groups $\tilde{H}_i(M \cap \breve{\Delta})$ (respectively, $I\tilde{H}_i(M \cap \breve{\Delta})$ with $i < \dim_{\mathbf{C}} M$) the isomorphic to them groups $\tilde{H}_i(M \cap \partial \bar{\Delta})$ or $\tilde{H}_{i+1}(M \cap \bar{\Delta}, M \cap \partial \bar{\Delta})$ (respectively, $I\tilde{H}_i(M \cap \Delta)$, see Proposition 1b) are substituted.* □

**Corollary 2:** *the last assertion of Theorem 1 about the formula (11).* □

## 4.2 Proof of Theorem 1

First we prove formula (9). Consider any smooth deformation contracting the disc $\Omega$ onto the segment $[-\sqrt{\delta}, \sqrt{\delta}]$. Using the fibre bundle (26), we can lift this deformation to the space $A \cap X \cap B^{(r+1)}$, which is therefore homotopy equivalent to the pre-image of this segment under the projection (26). Since the pre-images of the endpoints of this segment are contractible and over the interior points of the segment our projection is locally trivial, this pre-image is homotopy equivalent to the suspension of the fibre over any interior point, for instance of the fibre $\mathcal{F}_0 \equiv A \cap X \cap B^{(r)}$. This implies formula (9), and formula (7) follows from it by induction over $l = r - m + 1$.

To prove formula (10), consider the following filtration $\{\Psi_0 \subset \Psi_1 \subset \Psi_2\}$ of the disc $\Omega$ mod $\partial \Omega$. The term $\Psi_0$ consists of two points $\pm\sqrt{\delta}$, and $\Psi_1$ consists of two segments $[\sqrt{\delta}, \epsilon]$ and $[-\epsilon, -\sqrt{\delta}]$, so that the set $\Omega \setminus \Psi_1$ is a 2-cell; $\Psi_2 \equiv \Omega$. Using the projection (26) we lift this filtration onto the set $A \cap X \cap B^{(r+1)}$ and consider the spectral sequence generated by this filtration and calculating the relative (mod $\partial B^{(r+1)}$) homology of this set.

Table 1: The groups $E^1_{p,q}$ for the main (top part), antiinvariant (left-side) and invariant (right-side) spectral sequences

| $p$ | 0 | 1 | 2 | $\geq 3$ |
|---|---|---|---|---|
| $E^1_{p,q}$ | $(\check{H}_{q-1}(\partial\mathcal{F}))^2$ | $(\check{H}_q(\mathcal{F},\partial\mathcal{F}))^2$ | $\check{H}_q(\mathcal{F},\partial\mathcal{F})$ | 0 |

| $p$ | 0 | 1 | $\geq 2$ | 0 | 1 | 2 | $\geq 3$ |
|---|---|---|---|---|---|---|---|
| $E^1_{p,q}$ | $\check{H}_{q-1}(\partial\mathcal{F})$ | $\check{H}_q(\mathcal{F},\partial\mathcal{F})$ | 0 | $\check{H}_{q-1}(\partial\mathcal{F})$ | $\check{H}_q(\mathcal{F},\partial\mathcal{F})$ | $\check{H}_q(\mathcal{F},\partial\mathcal{F})$ | 0 |

**Proposition 7.** *For any $q$, the groups $E^1_{p,q}$ of this spectral sequence are as shown in the top part of Table 1 (where $\mathcal{F}$ is any generic fibre of (26), for instance $\mathcal{F}_0$).*

Indeed, the assertion about the column $\{p = 0\}$ follows from the fact that the fibres over the points $\pm\sqrt{\delta}$ are contractible, while their boundaries (= intersections with $\partial B^{(r+1)}$) are isotopic to these for nonsingular fibres. The assertions about columns $\{p = 1\}$ and $\{p = 2\}$ follow immediately from the constructions. □

From now on and up to § 4.5 we assume that the coefficient group allows division by 2.

The involution $\Omega \to \Omega$ sending $z$ to $-z$ acts on our spectral sequence and splits it into invariant ($E^{r+}_{p,q}$) and antiinvariant ($E^{r-}_{p,q}$) subsequences, which we can investigate separately.

All the column $\{p = 2\}$ of the term $E^1$ is invariant under this action, and the union of columns $\{p = 0\}$ and $\{p = 1\}$ splits into two parts, each of which coincides with the term $E^1$ of the spectral sequence calculating the homology mod $\partial B^{(r+1)}$ of the pre-image of either of the two components of the set $\Psi_1$. In particular, the terms $E^{1-}_{p,q}$ and $E^{1+}_{p,q}$ are as shown in the left and the right sides of the bottom row of Table 1, respectively.

The differential $\partial_1 : E^{1-}_{1,q} \to E^{1-}_{0,q}$ coincides with the boundary operator $\check{H}_q(\mathcal{F},\partial\mathcal{F}) \to \check{H}_{q-1}(\partial\mathcal{F})$, hence the antiinvariant spectral sequence is nothing but the exact sequence of the pair $(\mathcal{F},\partial\mathcal{F})$ calculating the absolute homology of $\mathcal{F}$. This gives the first term in the right-hand part of formula (10). It is easy to see that the realization of this term is exactly the one indicated in § 3.2.

The invariant spectral sequence is naturally isomorphic to the one calculating the homology group of the space $A \cap X \cap B^{(r+1)}$ reduced modulo the union of $\partial B^{(r+1)}$ and the half-space given by the condition $\mathrm{Re}\, z_{n-r} \leq -\delta$, in particular the final terms of both spectral sequences are adjoined to this homology.

The second summand in (10) (or, equivalently, in (8)) follows now from the next Proposition 8 (and our "invariant" spectral sequence is reduced to the exact sequence (29)).

**Proposition 8.** *Two quotient spaces $(A \cap X \cap B^{(r+1)})/(\partial B^{(r+1)} \cup \{z \in \mathbf{C}^n \mid \mathrm{Re}\, z_{n-r} \leq -\delta\})$ and $(A \cap B^{(r)})/\partial B^{(r)}$ are homotopy equivalent to one another.*

*Proof.* Consider the one-parametric family of $r$-dimensional affine subspaces $Y_\tau$, $\tau \in [0, \pi/2]$, in $\mathbf{C}^n$ containing the plane $X \cap T^r$ and distinguished by conditions $z_1 = \cdots = z_{n-r-1} = 0, \cos\tau \cdot (z_n - \delta) + \sin\tau \cdot z_{n-r} = 0$, so that $Y_0 = X \cap T^{r+1}$, $Y_{\pi/2} = T^r$. The family of quotient spaces $(A \cap Y_\tau \cap B^{(r+1)})/(\partial B^{(r+1)} \cup \{z \in \mathbf{C}^n \mid \operatorname{Re} z_{n-r} \leq -\delta\})$ realizes then a homotopy between two quotient spaces compared in the proposition. $\square$

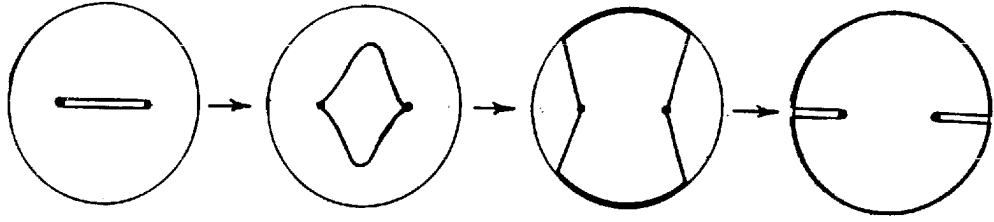

Figure 4: Deformations of the segment $[-\sqrt{\delta}, \sqrt{\delta}]$

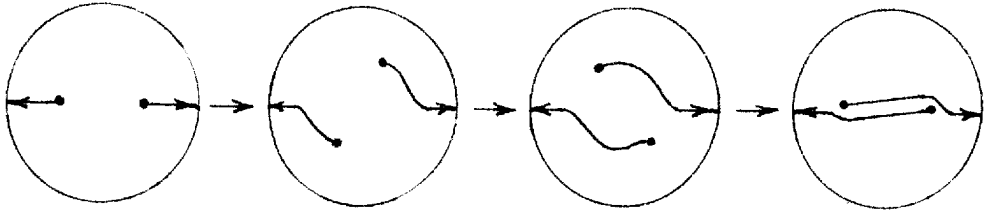

Figure 5: Variation of the disc $\Omega$

### 4.3 Proof of Theorem 2

The segment $[-\sqrt{\delta}, \sqrt{\delta}]$, considered as a relative cycle in $\Omega$ (mod $\Psi_0 \cup \partial\Omega$), can be deformed in two ways into the union of two segments $[-\sqrt{\delta}, -\epsilon]$ and $[\epsilon, \sqrt{\delta}]$; see Figure 4. Given a cycle $\alpha \in \tilde{H}_j(\mathcal{F})$, construct the cycle $\Sigma(\alpha) \in \mathcal{H}_{j+1}(r+1)$ as in § 3.2, divide it formally by 2, and deform these two halves in $A \cap X \cap B^{(r+1)}$ in correspondence with Figure 4: for any instant of this deformation these halves are realized by locally trivial bundles with standard fibres $\alpha/2$ over the points of the two curves obtained from the interval $(-\sqrt{\delta}, \sqrt{\delta})$. By construction, the resulting cycle is homologous to $\Sigma(\alpha)$ (mod $\partial B^{(r+1)}$), and is equal to $-(\rightleftharpoons(\alpha) + \overline{\rightleftharpoons(\alpha)})$, where $\overline{\rightleftharpoons(\alpha)}$ is obtained from $\alpha$ in the same way as $\rightleftharpoons(\alpha)$ but by means of transportation over a path complex conjugate to the S-shaped path in Figure 2a. By Lemma 3, $\overline{\rightleftharpoons(\alpha)} \simeq \rightleftharpoons(\alpha + \operatorname{Var}_r \circ J_r(\alpha))$ and Theorem 2 is proved. $\square$

### 4.4 Proof of Theorems 3 and 10

When $\xi$ moves along the circle $\delta \cdot e^{i\tau}, \tau \in [0, 2\pi]$, the ramification points $\pm\sqrt{\xi} \subset \Omega$ of the fibre bundle $z_{n-r} : A \cap X_\xi \cap B^{(r+1)} \to \Omega$ move as is shown in Figure 5. We move the segments $[-\epsilon, -\sqrt{\xi}], [\sqrt{\xi}, \epsilon]$ connecting the ramification points with "infinity" in such a way that they coincide with the parts of the real axis close to the boundary of $\Omega$. Let us also move the cycle realizing $\rightleftharpoons(\alpha)$ in such a way that at any instant of this deformation it forms a fibre bundle with standard fibre $\alpha/2$ over the union of two corresponding segments (except for their endpoints). At the final instant of the monodromy each of these two parts of $\rightleftharpoons(\alpha)$ will increase by the cycle $\Sigma(\alpha)/2$ swept by the formal half of the initial cycle $\alpha$ in transport over the added part of the resulting segment, so that $\operatorname{Var}_{r+1}(\rightleftharpoons(\alpha)) = \Sigma(\alpha)$.

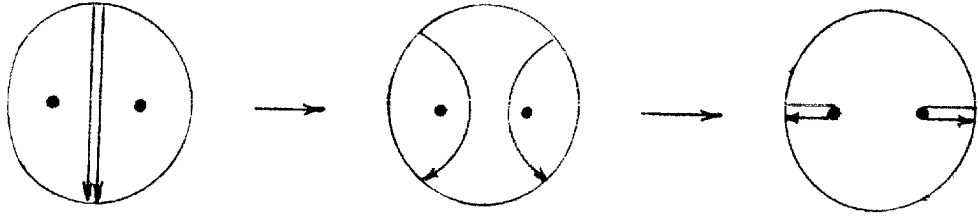

Figure 6: Homotopy between the cycles stab($\gamma$) and $\rightleftharpoons \circ \operatorname{Var}(\gamma)$

The restriction of $\operatorname{Var}_{r+1}$ to the second summand in (10) is a zero operator because all elements of this summand are invariant under the action of the involution $\{z_{n-r} \to -z_{n-r}\}$ and all elements of the target group $\mathcal{H}_*(r+1)$ are antiinvariant. Theorem 3 is proved. □

Theorem 10 is proved in the same way as Theorem 2, see Figure 6. □

## 4.5 On the proofs of Theorems 1–4, 10 in the case of general coefficients

In the case of general coefficients, we have not a splitting of the main spectral sequence into the invariant and antiinvariant parts: it splits only in a semidirect sum of the "invariant" subsequence and "antiinvariant" quotient sequence, which are isomorphic to the ones considered in § 4.2. The geometrical realization of the generators of the invariant sequence is the same as earlier, and the term $E^1$ of the antiinvariant one can be generated by the cycles lying only over the segment $[\sqrt{\delta}, \epsilon]$. Thus we have an exact sequence

$$\cdots \to \tilde{H}_{j+1}(A' \cap T^r, \partial B) \to \bar{\mathcal{H}}_{j+1}(r+1) \to \mathcal{H}_j(r) \to \cdots;$$

it splits because of the canonical realization of the generating cycles in $\bar{\mathcal{H}}_{j+1}(r+1)$ indicated in § 3.2.

Proposition 3 and Theorems 2, 3 and 10 follow immediately from this realization and Lemma 3. Indeed, if $\alpha = \operatorname{Var}_r(\gamma)$, $\gamma \in \bar{\mathcal{H}}_j(r)$, then the homology between the cycles $\rightleftharpoons(\alpha)$ and $\rightharpoonup(\alpha)$ is provided by the two-dimensional family of relative cycles in the fibres of the bundle (26) over the points of $\Omega \setminus \Psi_1$, obtained by the Gauss–Manin transportation from a cycle in $\mathcal{F}_0$ realizing $\gamma/2$. Theorem 2 (respectively, 3, respectively, 10) follows from the move of only the upper example of the segment $[-\sqrt{\delta}, \sqrt{\delta}]$ in Figure 4 (respectively, only the segment $[\sqrt{\delta}, \epsilon]$ in Figure 5, respectively, only the right-hand example of the imaginary axis in Figure 6). □

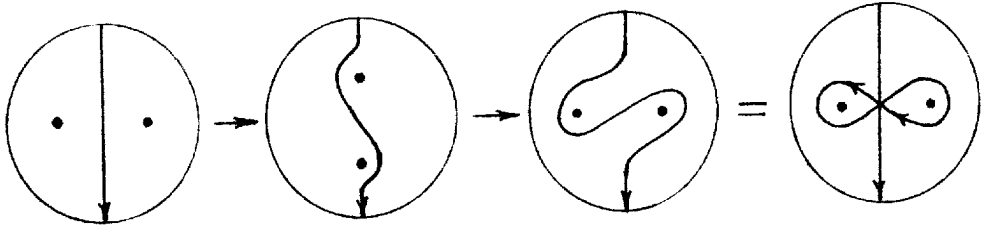

Figure 7: Stabilization of the variation operator in the intersection homology

## 5 Proofs for the intersection homology

**Theorem 13.** *For any $l = 1, \ldots, k$ and $\alpha \in I\bar{\mathcal{H}}_{d+l-1}(m+l-1)$ the homomorphism $\mathrm{Var}_{m+l}$ maps $\mathrm{stab}(\alpha)$ into $\infty(\alpha)$.*

*Proof:* see Figure 7.

**Corollary 1:** *Theorem 7.* □

**Corollary 2:** *formula (23) of Theorem 6.*

Indeed, if $l = 1$, then this formula coincides with Corollary 1 of Lemma 4. If $l > 1$, then $\langle \mathrm{stab}^l(\alpha), \infty^l(\beta) \rangle = \langle \mathrm{stab} \circ \mathrm{stab}^{l-1}(\alpha), \infty \circ \mathrm{stab}^{l-1}(\beta) \rangle = (-1)^{d+l-1} \langle \mathrm{stab}^{l-1}(\alpha), \mathrm{Var}_{d+l-1} \circ \mathrm{stab}^{l-1}(\beta) \rangle = (-1)^{d+l-1} \langle \mathrm{stab}^{l-1}(\alpha), \infty \circ \mathrm{stab}^{l-2}(\beta) \rangle = (-1)^{d+l-1} \langle \mathrm{stab}^{l-1}(\alpha), \infty^{l-1}(\beta) \rangle$; here the first and the last identities are the definition of $\infty^i$, the second one is Corollary 1 of Lemma 4, and the third is Theorem 13. □

**Corollary 3.** *For any natural $l \leq k$ the operator $\mathrm{Var}_{m+l} : I\bar{\mathcal{H}}_{d+l}(m+l) \to I\mathcal{H}_{d+l}(m+l)$ is an isomorphism.*

Indeed, by Corollary 2 of Theorem 12 these two groups are generated respectively by the elements of the form $\mathrm{stab}(\alpha)$ and $\infty(\alpha)$ and have equal dimensions. □

The formulae (15), (16) of Theorem 5 follow immediately from this corollary and Corollary 2 of Theorem 12.

**Theorem 14.** *For any $l = 1, \ldots, k$ and $\alpha \in I\bar{\mathcal{H}}_{d+l-1}(m+l-1)$ the obvious homomorphism*
$$J_{m+l} : I\mathcal{H}_{d+l}(m+l) \to I\bar{\mathcal{H}}_{d+l}(m+l)$$
*maps the element $\infty(\alpha)$ into $-\mathrm{stab}(2\alpha + J_{m+l-1} \circ \mathrm{Var}_{m+l-1}(\alpha))$.*

*Proof:* see Figure 8. □

Theorem 8 follows immediately from Theorems 13 and 14. □

**Theorem 15.** *Let $l$ be a natural number, $l \leq k$. Then for any $i \leq d+l-2$,*

$$I\mathcal{H}_i(m+l) \simeq I\mathcal{H}_i(m+l-1), \tag{31}$$

$$I\bar{\mathcal{H}}_{2d+2l-i}(m+l) \simeq I\bar{\mathcal{H}}_{2d+2l-i-2}(m+l-1); \tag{32}$$

$$I\mathcal{H}_{d+l-1}(m+l) \simeq I\mathcal{H}_{d+l-1}(m+l-1)/\mathrm{Im}\,\mathrm{Var}_{m+l-1}(I\bar{\mathcal{H}}_{d+l-1}(m+l-1)), \tag{33}$$

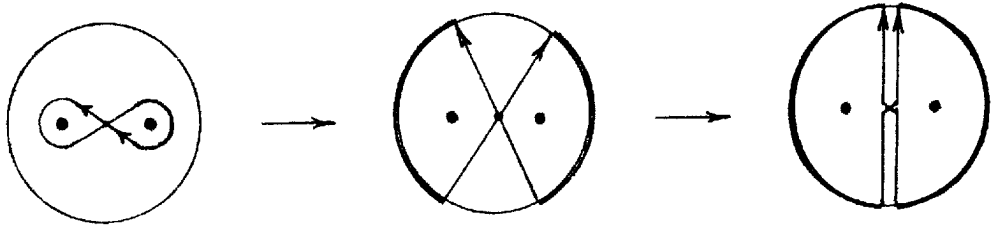

Figure 8: Stabilized map of the absolute intersection homology into the relative one

$$I\bar{\mathcal{H}}_{d+l+1}(m+l) \simeq \operatorname{Ker} \operatorname{Var}_{m+l-1}(I\bar{\mathcal{H}}_{d+l-1}(m+l-1)); \tag{34}$$

*the first and third of these isomorphisms are induced by the identical imbeddings, and the second and fourth are dual to them.*

The formulae (17)–(22) follow from these, from Proposition 2 and Corollary 3 of Theorem 13 by induction over $l$.

*Proof.* Let $W$ be the union of two very small (with respect to $\epsilon$ and $\delta$) open discs $W_+$ and $W_-$ in $\mathbf{C}^n$ centred at the distinguished points of the projection (26). Set $q = m + l$, $\Xi \equiv A' \cap X \cap B^{(q)}$, and consider the exact sequence of the pair $(\Xi, \Xi \cap W)$,

$$\cdots \to IH_{i+1}(\Xi, \Xi \cap W) \to IH_i(\Xi \cap W) \to I\mathcal{H}_i(q) \to \cdots, \tag{35}$$

and the Poincaré dual to it sequence of the triple $(\Xi, \Xi \setminus W, \Xi \cap \partial B')$,

$$\cdots \leftarrow IH_{2d+2l-i-1}((\Xi \setminus W), \partial B') \leftarrow IH_{2d+2l-i}(\Xi, \Xi \setminus W) \leftarrow I\bar{\mathcal{H}}_{2d+2l-i}(q) \leftarrow \cdots. \tag{36}$$

For $i \leq d + l - 2$ the marginal arrows in (36) (and hence also in (35)) are trivial. Indeed, arguing as in the proof of Proposition 6 or Theorem 11, we obtain that any element of the left group in (36) is equal to $\operatorname{stab}(\gamma)$ for some $\gamma \in I\bar{\mathcal{H}}_{d+2l-i-2}(q-1)$ and hence, by Theorem 12, is homologous to zero. Let us study the obtained short exact sequence (35).

The group $IH_{i+1}(\Xi, \Xi \cap W)$ for such $i$ is obviously isomorphic to $IH_i(\mathcal{F}_0)$ (the isomorphism is realized by the suspension operation like in the construction of cycles $\Sigma(\alpha)$, see also the proof of Theorem 1). The group $IH_i(\Xi \cap W)$ is isomorphic to $(IH_i(\mathcal{F}_0))^2$, since in the exact sequence (30) for the homology of spaces $\Xi \cap W_+$ and $\Xi \cap W_-$ (applicable by Corollary 1 of Lemma 5) the third term (and the term "before the first arrow") is trivial by dimensional reason and the first term equals $IH_i(\mathcal{F}_0)$. The second arrow in (35) connecting these groups is thus isomorphic to the diagonal imbedding $IH_i(\mathcal{F}_0) \to (IH_i(\mathcal{F}_0))^2$ and formula (31) follows; (32) is just a dual of it.

If $i = d + l - 1$, the nonsplitting fragment of the sequence (35) is one term longer:

$$0 \to I\mathcal{H}_{d+l}(q) \to IH_{d+l}(\Xi, \Xi \cap W) \xrightarrow{\partial} IH_{d+l-1}(\Xi \cap W) \to I\mathcal{H}_{d+l-1}(q) \to 0, \tag{37}$$

where the right-hand zero is proved as above and the left zero follows from Proposition 1b. By Corollary 2 of Theorem 12 the map $I\mathcal{H}_{d+l}(q) \to IH_{d+l}(\Xi, \Xi \cap W)$ in this sequence is isomorphic to the inclusion Im $\text{Var}_{q-1}(I\bar{\mathcal{H}}_{d+l-1}(q-1)) \hookrightarrow I\mathcal{H}_{d+l-1}(q-1)$, and by (30) the term $I\mathcal{H}_{d+l-1}(\Xi \cap W)$ is isomorphic to $[I\bar{\mathcal{H}}_{d+l-1}(q-1)/\text{Im Var}_{q-1}(I\bar{\mathcal{H}}_{d+l-1}(q-1))]^2$; by the construction of all these isomorphisms, the image of the map $\partial$ in (37) is the diagonal in this square group. Theorem 15 is completely proved. $\square$

Two first identities of Theorem 9 and the third identity for $l = 1$ follow immediately from the realization of operators $\Sigma$, stab and $\infty$; for arbitrary $l$ the last identity follows by induction over $l$ from the case $l = 1$ and Theorem 13. $\square$

**Acknowledgement.** I thank V. A. Ginzburg for consultations on intersection homology theory.